\documentclass{article}
\usepackage{amsmath}

\newtheorem{theorem}{Theorem}
\newtheorem{acknowledgement}[theorem]{Acknowledgement}

\newtheorem{corollary}[theorem]{Corollary}

\newtheorem{example}[theorem]{Example}

\newtheorem{lemma}[theorem]{Lemma}

\newtheorem{remark}[theorem]{Remark}

\newenvironment{proof}[1][Proof]{\noindent\textbf{#1.} }{\ \rule{0.5em}{0.5em}}
\include {mak}
\input{tcilatex}

\begin{document}

\title{The role of binomial type sequences in determination identities for
Bell polynomials}
\author{Miloud Mihoubi \\
\ \ \ \ \ \ \ \ \ \ \ \ \ \ \ \ \\
U.S.T.H.B., Faculty of Mathematics, \\
B.P. 32, El-Alia, 16111, Bab Ezzouar,\\
Algiers, Algeria.}
\maketitle

\begin{quote}
\textbf{Abstract. }Our paper deals about identities involving Bell
polynomials. Some identities on Bell polynomials derived using generating
function and successive derivatives of binomial type sequences. We give some
relations between Bell polynomials and binomial type sequences in first
part, and, we generalize the results obtained in [4] in second part.

\ \ \ \ \ \ \ \ \ \ \ \ \ \ \ \ \ \ \ \ 

\end{quote}

\section{Introduction}

Recall that the (exponential) partial Bell polynomials $B_{n,k}\left(
x_{1},x_{2},..\right) $\ are defined by their generating function:%
\begin{equation}
\underset{n=k}{\overset{\infty }{\sum }}B_{n,k}\left( x_{1},x_{2},...\right) 
\frac{t^{n}}{n!}=\frac{1}{k!}\left( \underset{m=1}{\overset{\infty }{\sum }}%
x_{m}\frac{t^{m}}{m!}\right) ^{k}.  \label{a}
\end{equation}%
and the (exponential) complete Bell polynomials $A_{n}\left(
x_{1},x_{2},..\right) $\ are given by%
\begin{equation*}
A_{n}\left( x_{1},x_{2},..\right) :=\underset{k=1}{\overset{n}{\sum }}%
B_{n,k}\left( x_{1},x_{2},...\right) \text{ \ \ with \ \ }A_{0}\left(
x_{1},x_{2},..\right) :=1.
\end{equation*}%
Comtet $\left[ 3\right] $ studies the Bell polynomials and gives some basic
properties for them. Some applications of Bell polynomials are given by
Riordan $\left[ 5\right] $ in combinatorial analysis and by S. Roman $\left[
6\right] $ in umbral calculus. Recently, by using the Lagrange inversion
formula (LIF), Abbas and Bouroubi $\left[ 1\right] $ give some identities
for the partial Bell polynomials, and, Mihoubi $\left[ 4\right] $ also gives
some extensions involving to the partial and complete Bell polynomials. For
any sequence $\left( x_{n};n\geq 1\right) $\ with $x_{1}=1$ and any natural
numbers $r,s,$ recall that Proposition 4 in $\left[ 4\right] $ gives:%
\begin{gather}
B_{n,k}\left( B\left( s,s\right) ,...,\frac{ms}{r\left( m-1\right) +s}\frac{%
B\left( \left( r+1\right) \left( m-1\right) +s,\;r\left( m-1\right)
+s\right) }{\binom{\left( r+1\right) \left( m-1\right) +s}{r\left(
m-1\right) +s}},...\right)  \label{w1} \\
=\binom{n}{k}\frac{sk}{r\left( n-k\right) +sk}\frac{B\left( \left(
r+1\right) \left( n-k\right) +sk,\;r\left( n-k\right) +sk\right) }{\binom{%
\left( r+1\right) \left( n-k\right) +sk}{r\left( n-k\right) +sk}},  \notag \\
\ \text{where }B\left( n,k\right) :=B_{n,k}\left(
x_{1},x_{2},x_{3},...\right) .  \notag
\end{gather}%
Then if we put%
\begin{equation}
Y\left( n,k\right) :=\binom{n}{k}\frac{sk}{r\left( n-k\right) +sk}\frac{%
B\left( \left( r+1\right) \left( n-k\right) +sk,\;r\left( n-k\right)
+sk\right) }{\binom{\left( r+1\right) \left( n-k\right) +sk}{r\left(
n-k\right) +sk}},  \label{w2}
\end{equation}%
we conclude, from (\ref{w1}), that the sequence $\left( Y\left( n,k\right)
\right) $ satisfies the equation:%
\begin{equation}
B_{n,k}\left( Y\left( 1,1\right) ,Y\left( 2,1\right) ,Y\left( 3,1\right)
,...\right) =Y\left( n,k\right) .  \label{h}
\end{equation}%
In general, from the definition (\ref{a}), if we put $\psi \left( t\right) =%
\underset{m=1}{\overset{\infty }{\sum }}x_{m}\dfrac{t^{m}}{m!}$\ and $%
Y\left( n,k\right) =\frac{1}{k!}D_{t=0}^{n}\left( \psi \left( t\right)
\right) ^{k},$\ then for all integers $n,k,$\ with $n\geq k\geq 1,$\ the
sequence $\left( Y\left( n,k\right) \right) $ satisfies (\ref{h}). Hence, to
find identities for the partial Bell polynomials, it suffices to find
sequences $\left( Y\left( n,k\right) \right) $\ satisfy the equation (\ref{h}%
).\newline
Similarly to the partial Bell polynomials, for any sequence $\left(
x_{n};n\geq 1\right) $\ with $x_{1}=1$ and any natural numbers $r,s$ $\left(
r\geq 1\right) ,$ another relation for the complete Bell polynomials is
given by Proposition 8 in $\left[ 4\right] $ by:%
\begin{eqnarray}
A_{n}\left( \frac{sB\left( r+1,\ r\right) }{r\left( r+1\right) },...,\frac{s%
}{nr}\frac{B\left( \left( r+1\right) n,\ nr\right) }{\binom{\left(
r+1\right) n}{nr}}\right) &=&\frac{s}{nr+s}\frac{B\left( \left( r+1\right)
n+s,\ nr+s\right) }{\binom{\left( r+1\right) n+s}{nr+s}}  \label{w3} \\
\text{where }B\left( n,k\right) &:&=B_{n,k}\left(
x_{1},x_{2},x_{3},...\right) ,  \notag
\end{eqnarray}%
and if we put 
\begin{equation}
Z\left( n,s\right) :=\frac{1}{nr+s}\frac{B\left( \left( r+1\right) n+s,\
nr+s\right) }{\binom{\left( r+1\right) n+s}{nr+s}},  \label{w4}
\end{equation}%
we conclude, from (\ref{w3}), that the sequence $\left( Z\left( n,s\right)
\right) $ satisfies the equation:%
\begin{equation}
A_{n}\left( sZ\left( 1,0\right) ,sZ\left( 2,0\right) ,...,sZ\left(
n,0\right) \right) =sZ\left( n,s\right) .  \label{b}
\end{equation}%
Hence, to find identities for the complete Bell polynomials, it suffices to
find sequences $\left( Z\left( n,s\right) \right) $\ satisfy the equation (%
\ref{b}). Therefore, to determine solutions for (\ref{h}) and (\ref{b}), we
exploit the strong connection between Bell polynomials and binomial type
sequences. For any binomial type sequence $\left( f_{n}\left( x\right)
\right) ,$ with $f_{0}\left( x\right) :=1,$ one of such connections is given
in $\left[ 6,\ p.\ 82\right] $\ by%
\begin{equation}
f_{n}\left( x\right) =\underset{k=1}{\overset{n}{\sum }}B_{n,k}\left(
x_{1},x_{2},...\right) x^{k}\text{ \ \ with \ }x_{n}=\frac{d}{dx}\left.
f_{n}\left( x\right) \right\vert _{x=0}.  \label{f}
\end{equation}%
On the basis of the results obtained in $\left[ 4\right] $ and the relation (%
\ref{f}), we derive in this paper some interesting identities and relations
related Bell polynomials and binomial type sequences. \newline
For the next of this paper, we will denote by $D_{t}f\left( t\right) ,\
D_{t}^{j}f\left( t\right) ,$\ $D_{t=x}f\left( t\right) $\ or $Df\left(
x\right) $ and $D_{t=x}^{j}f\left( t\right) ,$\ respectively, for the
derivative of $f,\ $the $j-th$\ derivative of $f,$\ the derivative of \ $f$\
evaluated at $t=x$, and the $j-th$\ derivative of\ $f$\ evaluated at $t=x.$

\section{Main results}

For the next of this work, for any sequence $\left( f_{n}\left( x\right)
\right) $\ of binomial type with $f_{0}\left( x\right) =1$\ and for a given
real number $a,$\ we define%
\begin{equation}
f_{n}\left( x;a\right) :=\frac{x}{an+x}f_{n}\left( an+x\right) \text{ \ with%
\textrm{\ \ }}f_{0}\left( x;a\right) =1.  \label{x0}
\end{equation}%
The sequence $\left( f_{n}\left( x;a\right) \right) $\ is of binomial type,
see $\left[ 4\right] .$\newline
To simplify any expression below we put%
\begin{gather*}
T=T\left( n,k\right) :=r\left( n-k\right) +sk,\text{ \ \ }R=R\left(
n,s\right) :=nr+s \\
\text{and }\left( f_{n}\left( x\right) \right) \text{ denotes a sequence of
binomial type with }f_{0}\left( x\right) =1.
\end{gather*}%
The two following theorems give some interesting relations between Bell
polynomials and binomial type sequences. These relations are used to deduce
several identities for partial and complete Bell polynomials as it is
illustrated below. To prove these theorems, we use the following Lemma:

\begin{lemma}
Let $n,k$\ be integers with $n\geq k\geq 1$ and $a,\alpha $ be a real
numbers. We have%
\begin{gather}
B_{n,k}\left( \alpha ,2D_{z=0}\left( e^{\alpha z}f_{1}\left( x+z;a\right)
\right) ,...,mD_{z=0}\left( e^{\alpha z}f_{m-1}\left( kx+z;a\right) \right)
,...\right)  \label{alpha} \\
=\binom{n}{k}D_{z=0}^{k}\left( e^{\alpha z}f_{n-k}\left( kx+z;a\right)
\right) .  \notag
\end{gather}%
This identity can be replaced when $\alpha =0$ by%
\begin{equation}
B_{n,k}\left( D_{x}f_{1}\left( x;a\right) ,...,D_{x}f_{m}\left( x;a\right)
,...\right) =\frac{1}{k!}D_{z=0}^{k}f_{n}\left( kx+z;a\right) .
\label{alpha1}
\end{equation}
\end{lemma}

\begin{theorem}
\label{1}Let $a,x,\alpha $\ be real numbers and $n,k,r,s$\ be natural
numbers with$\ n\geq k\geq 1$ and $r+s\geq 1.$\ Then the sequence%
\begin{equation}
Y\left( n,k\right) :=\binom{n}{k}\frac{sk}{T}D_{z=0}^{T}\left( e^{\alpha
z}f_{n-k}\left( Tx+z;a\right) \right)  \label{x11}
\end{equation}%
satisfies (\ref{h}). For $\alpha =0$ the sequence $\left( Y\left( n,k\right)
\right) $ can be replaced by%
\begin{equation}
Y\left( n,k\right) :=\frac{n!}{k!\left( T+n-k\right) !}\frac{sk}{T}%
D_{z=0}^{T}f_{T+n-k}\left( Tx+z;a\right) .  \label{x1}
\end{equation}%
For $r=s=0,$ we put $Y\left( n,k\right) :=\binom{n}{k}f_{n-k}\left(
x;a\right) .$
\end{theorem}

\begin{theorem}
\label{2}Let $a,x,\alpha $\ be real numbers and $n,r,s$\ be natural numbers
with$\ n\geq 1,$\ $r\geq 1.$\ Then for $\alpha \neq 0$ the sequence%
\begin{equation}
Z\left( n,s\right) :=\frac{1}{\alpha ^{s}}\frac{1}{R}D_{z=0}^{R}\left(
e^{\alpha z}f_{n}\left( Rx+z;a\right) \right)  \label{x22}
\end{equation}%
satisfies (\ref{b}), and (case $\alpha =0$) the sequence%
\begin{equation}
Z\left( n,s\right) :=\frac{1}{\left( Df_{1}\left( 0\right) \right) ^{s}}%
\frac{1}{R}\frac{n!}{\left( R+n\right) !}D_{z=0}^{R}f_{R+n}\left(
Rx+z;a\right)  \label{x2}
\end{equation}%
satisfies (\ref{b}).
\end{theorem}

More generally, Theorem \ref{1} can be generalized as follows:

\begin{theorem}
\label{3}Let $\left( a_{n};n\geq 1\right) $\ be a real sequence; $%
n,k,r,s,u,v $\ be natural numbers with\ $n\geq k\geq 1,$\ $r+s\geq 1$\ and $%
x,a,\alpha ,\beta ,\lambda $\ be real numbers. Then\ the sequence%
\begin{equation}
Y\left( n,k\right) :=\binom{n}{k}\frac{sk}{T}T!\underset{j\geq T}{\sum }%
B_{j,T}\left( a_{1},a_{2},...\right) D_{z=\beta j+\lambda T}^{ju+vT}\left\{
e^{\alpha z}f_{n-k}\left( z;a\right) \right\} \frac{x^{j}}{j!}  \label{x55}
\end{equation}%
satisfies (\ref{h}). For $\alpha =0$ the above sequence can be replaced by%
\begin{equation}
Y\left( n,k\right) :=\frac{n!}{k!}\frac{sk}{T}T!\underset{j=T}{\overset{h+T}{%
\sum }}B_{j,T}\left( a_{1},a_{2},...\right) \frac{D_{z=\beta j+\lambda
T}^{ju+vT}f_{\left( u+v\right) T+n-k}\left( z;a\right) }{\left( \left(
u+v\right) T+n-k\right) !}\frac{x^{j}}{j!}  \label{x5}
\end{equation}%
where $h=\left[ \frac{n-k}{u}\right] $\ for $u\geq 1,$\ $h=\infty $\ for $%
u=0 $\ and $\left[ x\right] $\ is the largest integer $\leq x.$
\end{theorem}

More generally, Theorem \ref{2} can be generalized as follows:

\begin{theorem}
\label{4}Let $\left( a_{n};n\geq 1\right) $\ be a real sequence, $n,r,s,u,v$%
\ be natural numbers with\ $r\geq 1$\ and $a,\alpha ,\beta ,\lambda $\ be
real numbers. Then\ for $\alpha \neq 0,$ the sequence%
\begin{gather}
Z\left( n,s\right) :=\frac{R!}{\gamma ^{s}R}\underset{j\geq R}{\sum }%
B_{j,R}\left( a_{1},a_{2},...\right) D_{z=\beta j+\lambda R}^{ju+vR}\left\{
e^{\alpha z}f_{n}\left( z;a\right) \right\} \frac{x^{j}}{j!}  \label{x77} \\
\text{with }\gamma :=\alpha ^{v}\varphi \left( x\alpha ^{u}\right) \text{ \
and \ }\varphi \left( x\right) :=\underset{i=1}{\overset{\infty }{\sum }}%
a_{i}\frac{x^{i}}{i!}  \notag
\end{gather}%
satisfies (\ref{b}). For $\alpha =0$ the above sequence can be replaced by%
\begin{gather}
Z\left( n,s\right) :=\frac{n!R!}{\gamma ^{s}R}\underset{j=R}{\overset{g+R}{%
\sum }}B_{j,R}\left( a_{1},a_{2},...\right) \frac{D_{z=\beta j+\lambda
R}^{ju+vR}f_{n+\left( u+v\right) R}\left( z;a\right) }{\left( n+\left(
u+v\right) R\right) !}\frac{x^{j}}{j!}  \label{x7} \\
\text{with }\gamma :=\left\{ 
\begin{array}{cc}
a_{1}x\left( Df_{1}\left( 0\right) \right) ^{u+v} & \text{if }u\geq 1 \\ 
\left( Df_{1}\left( 0\right) \right) ^{v}\varphi \left( x\right) & \text{if }%
u=0%
\end{array}%
\right. \text{\ and }\varphi \left( x\right) :=\underset{i=1}{\overset{%
\infty }{\sum }}a_{i}\frac{x^{i}}{i!},  \notag
\end{gather}%
where $g:=\left[ \frac{n}{u}\right] $ for $u\geq 1$, and, $g=\infty $ for $%
u=0.$
\end{theorem}

\begin{remark}
For $a_{n}=0$\ $\left( n\geq 2\right) $\ in Theorem \ref{3}\ we obtain
Theorem \ref{1}.\newline
For $\ a_{n}=0$\ $\left( n\geq 2\right) $\ \ in Theorem \ref{4}\ we obtain
Theorem \ref{2}.\newline
For $x=\alpha =0$\ in Theorem \ref{1}\ we obtain Corollary 5 in $\left[ 4%
\right] $.\newline
For $x=\alpha =0$\ in Theorem \ref{2} we get Corollary 9 in $\left[ 4\right]
.$\newline
From (\ref{x11}) the sequence $Y_{1}\left( n,k\right) :=\underset{\alpha
\rightarrow 0}{\lim }\alpha ^{-T}Y\left( n,k\right) $ satisfies (\ref{h})
and gives Proposition 1 in $\left[ 4\right] ,$ and similarly, from (\ref{x22}%
) the sequence $Z_{1}\left( n,s\right) :=\underset{\alpha \rightarrow 0}{%
\lim }\alpha ^{-nr}Z\left( n,s\right) $ satisfies (\ref{b}) and gives
Proposition 3 in $\left[ 4\right] .$\newline
By using (\ref{f}) and (\ref{x0}) we can construct several binomial type
polynomials as 
\begin{equation*}
p_{n}\left( t\right) :=t\overset{n}{\underset{k=1}{\sum }}Y\left( n,k\right)
\left( bn+t\right) ^{k-1},\ \left( p_{0}\left( t\right) :=1\right) ,
\end{equation*}%
where $\left( Y\left( n,k\right) \right) $\ is given by (\ref{x1}) or (\ref%
{x11}).
\end{remark}

\section{Applications}

We give in this section another versions of Theorems \ref{1}, \ref{2} and we
present some particular cases of the above results.

\subsection{Some\ applications\ of\ Theorem \protect\ref{1}}

The following corollaries gives a practical version of Theorem \ref{1}.

\begin{corollary}
Under the hypothesis of Theorem \ref{1} the sequence%
\begin{equation}
Y\left( n,k\right) :=\binom{n}{k}\frac{sk}{T}\underset{j=0}{\overset{n-k}{%
\sum }}\binom{T}{j}D_{z=Tx}^{j}\ f_{n-k}\left( z;a\right) \alpha ^{j}
\label{p1}
\end{equation}%
satisfies (\ref{h}).
\end{corollary}

\begin{proof}
From (\ref{x11}), we have%
\begin{eqnarray*}
Y\left( n,k\right) &=&\binom{n}{k}\frac{sk}{T}\underset{j=0}{\overset{n-k}{%
\sum }}\binom{T}{j}D_{z=0}^{T-j}\left( e^{z/\alpha }\right)
D_{z=0}^{j}f_{n-k}\left( Tx+z;a\right) \\
&=&\binom{n}{k}\frac{sk}{T}\underset{j=0}{\overset{n-k}{\sum }}\binom{T}{j}%
D_{z=0}^{j}f_{n-k}\left( Tx+z;a\right) \left( \frac{1}{\alpha }\right)
^{T-j}.
\end{eqnarray*}%
We conclude noticing that $Y_{1}\left( n,k\right) :=\alpha ^{T}Y\left(
n,k\right) $ satisfies (\ref{h}).
\end{proof}

\begin{example}
For $f_{n}\left( x\right) =x^{n}$ the sequence $\left( Y\left( n,k\right)
\right) $ in (\ref{p1}) becomes%
\begin{equation*}
Y\left( n,k\right) :=\binom{n}{k}\frac{sk}{T}\underset{j=0}{\overset{n-k}{%
\sum }}\binom{T}{j}\frac{\left( n-k\right) !}{\left( n-k-j\right) !}\left(
Tx+a\left( n-k\right) \right) ^{n-k-j-1}\left( Tx+aj\right) \alpha ^{j}
\end{equation*}%
and for $a=0,$ $\alpha =1$ the last sequence becomes 
\begin{equation*}
Y\left( n,k\right) :=T!\binom{n}{k}\frac{sk}{T}\underset{j=0}{\overset{n-k}{%
\sum }}\binom{n-k}{j}\frac{\left( Tx\right) ^{n-k-j}}{\left( T-j\right) !}.
\end{equation*}
\end{example}

The following corollary gives a practical version of Theorem \ref{1} when $%
\alpha =0$.

\begin{corollary}
Under the hypothesis of Theorem \ref{1} the sequence%
\begin{gather}
Y\left( n,k\right) :=\frac{sk}{T}\binom{n}{k}\binom{T+n-k}{n-k}^{-1}\underset%
{j=0}{\overset{n-k}{\sum }}\binom{T+j-1}{T-1}B\left( T+n-k,T+j\right) \times
\label{z1} \\
\left( \left( \left( r+1\right) c-b\right) j+b\left( n-k\right) +csk\right)
\left( b\left( n-k\right) +csk\right) ^{T+j-1}  \notag
\end{gather}%
satisfies (\ref{h}), and in particular the sequence%
\begin{eqnarray}
Y\left( n,k\right) &:&=\frac{sk}{T}\binom{n}{k}\binom{T+n-k}{n-k}^{-1}\times
\label{z2} \\
&&\underset{j=0}{\overset{n-k}{\sum }}\binom{T+j-1}{T-1}B\left(
T+n-k,T+j\right) \left( T+n-k\right) ^{j}x^{j}  \notag
\end{eqnarray}%
satisfies (\ref{h}), where $B_{n,k}\left( x_{1},x_{2},x_{3},..\right)
:=B\left( n,k\right) .$
\end{corollary}

\begin{proof}
Let $\left( x_{n};n\geq 1\right) $\ be a sequence of real numbers. From (\ref%
{x1}), it suffices to express $Y\left( n,k\right) $\ by considering the
binomial type sequence 
\begin{equation}
f_{n}\left( x;a\right) :=x\underset{k=1}{\overset{n}{\sum }}B_{n,k}\left(
x_{1},x_{2},...\right) \left( an+x\right) ^{k-1}\ \text{with }f_{0}\left(
x\right) =1,  \label{p5}
\end{equation}%
put after $b:=ra+rx+a,\ \ c:=x+a.$ To obtain (\ref{z2}), it suffices to
choice $b=\left( r+1\right) c,\ c=x$\ in (\ref{z1}). To obtain (\ref{z2}),
it suffices to choice $b=\left( r+1\right) c,\ c=x$\ in (\ref{z1}).
\end{proof}

\begin{example}
By using the well-known identity $B_{n,k}\left( 1,2,3,...\right) =\binom{n}{k%
}k^{n-k},$\ if $x_{n}=n,$\ the sequence given by (\ref{z2}) becomes 
\begin{equation*}
Y\left( n,k\right) =\frac{n!}{k!}\frac{sk}{\left( T+n-k\right) !}\underset{%
j=0}{\overset{n-k}{\sum }}\binom{T+n-k}{T+j}\left( T+j\right) ^{n-k-j}\left(
T+n-k\right) ^{j}\frac{x^{j}}{j!},
\end{equation*}%
and by using the well-known identity $B_{n,k}\left( 1!,2!,3!,...\right) =%
\binom{n}{k}\frac{\left( n-1\right) !}{\left( k-1\right) !}$\ (the Lah
numbers), if $x_{n}=n!,$\ the sequence given by (\ref{z2}) becomes%
\begin{equation*}
Y\left( n,k\right) =\frac{n!}{k!}\frac{sk}{T+n-k}\underset{j=0}{\overset{n-k}%
{\sum }}\binom{T+n-k}{T+j}\left( T+n-k\right) ^{j}\frac{x^{j}}{j!}.
\end{equation*}%
More applications can be constructed by choosing in (\ref{z1}) or (\ref{z2}) 
$x_{n}$ as $1,$\ $nf_{n-1}\left( y;a\right) ,$\ $D_{x}f_{n}\left( x;a\right)
.$
\end{example}

\begin{example}
For $r=0,$\ $s=1,\ \alpha =0$\ and $f_{n}\left( x\right) :=B_{n}\left(
x\right) =\underset{j=1}{\overset{n}{\sum }}\binom{n}{j}S\left( n,j\right) $%
\ in (\ref{x1}), with $S\left( n,j\right) $\ are the Stirling numbers of the
second kind, we obtain:%
\begin{gather*}
B_{n,k}\left( 1,...,\frac{xB_{m+1}\left( am+x\right) -\left(
x^{2}+amx-am\right) B_{m}\left( am+x\right) }{\left( am+x\right) ^{2}}%
,...\right) \\
=kx\underset{l=0}{\overset{n-k}{\sum }}\binom{l+k}{k}S\left( n,l+k\right)
\left( an+kx\right) ^{l-1}.
\end{gather*}%
For $r=0,$\ $s=1,\ \alpha =0$\ and \ $f_{n}\left( x\right) :=\left[ x\right]
_{n}=\underset{j=1}{\overset{n}{\sum }}s\left( n,j\right) x^{j}$\ in (\ref%
{x1}), where $s\left( n,j\right) $\ are the Stirling numbers of the first
kind, we obtain:%
\begin{equation*}
\left. 
\begin{array}{c}
B_{n,k}\left( 1,...,\left[ am+x-1\right] _{m-1}\left( m-\underset{i=1}{%
\overset{m-1}{\sum }}\frac{am-i}{am+x-i}\right) ,...\right) = \\ 
kx\underset{l=0}{\overset{n-k}{\sum }}\binom{l+k}{k}s\left( n,l+k\right)
\left( an+kx\right) ^{l-1}.%
\end{array}%
\right.
\end{equation*}%
For $r=0,$\ $s=1,\ \alpha =0$\ and \ $f_{n}\left( x\right) :=\left[ x\right]
^{n}=\underset{j=1}{\overset{n}{\sum }}\left\vert s\left( n,j\right)
\right\vert x^{j}$\ in (\ref{x1}), where $\left\vert s\left( n,j\right)
\right\vert $\ are the absolute Stirling numbers of the first kind, we
obtain:%
\begin{gather*}
B_{n,k}\left( 1,...,\left[ am+x-1\right] ^{m-1}\left( m-\underset{i=1}{%
\overset{m-1}{\sum }}\frac{am+i}{am+i+x}\right) ,...\right) \\
=kx\underset{l=0}{\overset{n-k}{\sum }}\binom{l+k}{k}\left\vert s\left(
n,l+k\right) \right\vert \left( an+kx\right) ^{l-1}.
\end{gather*}
\end{example}

\subsection{Some\ applications\ of\ Theorem\ \protect\ref{2}}

The following corollaries gives a practical version of Theorem \ref{2}

\begin{corollary}
Under the hypothesis of Theorem \ref{2} the sequence%
\begin{equation}
Z\left( n,s\right) :=\frac{1}{R}\underset{j=0}{\overset{n}{\sum }}\binom{R}{j%
}D_{z=Rx}^{j}\ f_{n}\left( z;a\right) \alpha ^{j}  \label{p6}
\end{equation}%
satisfies (\ref{b}).
\end{corollary}

\begin{proof}
From (\ref{x22}) we have%
\begin{eqnarray*}
Z\left( n,s\right) &=&\frac{\alpha ^{s}}{R}D_{z=0}^{R}\left( e^{z/\alpha
}f_{n}\left( Rx+z;a\right) \right) \\
&=&\frac{\alpha ^{s}}{R}\underset{j=0}{\overset{n}{\sum }}\binom{R}{j}%
D_{z=0}^{j}\left( f_{n}\left( Rx+z;a\right) \right) \frac{1}{\alpha ^{R-j}}
\\
&=&\frac{1}{R}\underset{j=0}{\overset{n}{\sum }}\binom{R}{j}%
D_{z=0}^{j}\left( f_{n}\left( Rx+z;a\right) \right) \frac{\alpha ^{j}}{%
\alpha ^{nr}}.
\end{eqnarray*}%
\newline
We conclude noticing that $Z_{1}\left( n,s\right) :=\alpha ^{nr}Y\left(
n,k\right) $ satisfies (\ref{h}).
\end{proof}

\begin{example}
For $f_{n}\left( x\right) =x^{n}$ and $r\geq 1$ the sequence given by (\ref%
{p6}) becomes%
\begin{equation*}
Z\left( n,s\right) =\frac{1}{R}\underset{j=0}{\overset{n}{\sum }}\binom{R}{j}%
\frac{n!}{\left( n-j\right) !}\left( Rx+an\right) ^{n-j-1}\left(
Rx+aj\right) \alpha ^{j}
\end{equation*}%
and for $a=0,$ $x=1,$ the last sequence becomes%
\begin{equation*}
Z\left( n,s\right) =\frac{1}{R}\underset{j=0}{\overset{n}{\sum }}\binom{R}{j}%
\frac{n!}{\left( n-j\right) !}R^{n-j}\alpha ^{j}.
\end{equation*}
\end{example}

\begin{corollary}
Under the hypothesis of Theorem \ref{2}\ the sequence%
\begin{gather}
Z\left( n,s\right) :=\frac{1}{\left( Df_{1}\left( 0\right) \right) ^{s}}%
\frac{1}{R}\binom{R+n}{n}^{-1}\underset{j=0}{\overset{n}{\sum }}\binom{R+j-1%
}{R-1}B\left( R+n,R+j\right)  \label{z3} \\
\times \left( cn+bs\right) ^{j-1}\left\{ \left( \left( r+1\right) b-c\right)
j+cn+bs\right\}  \notag
\end{gather}%
satisfies (\ref{b}), and in particular the sequence%
\begin{eqnarray}
Z\left( n,s\right) &:&=\frac{1}{\left( Df_{1}\left( 0\right) \right) ^{s}}%
\frac{1}{R}\binom{R+n}{n}^{-1}\underset{j=0}{\overset{n}{\sum }}\binom{R+j-1%
}{R-1}\times  \label{z4} \\
&&B\left( R+n,R+j\right) \left( R+n\right) ^{j}x^{j}  \notag
\end{eqnarray}%
satisfies (\ref{b}), where $B_{n,k}\left( x_{1},x_{2},x_{3},..\right)
:=B\left( n,k\right) .$
\end{corollary}

\begin{proof}
Let $\left( x_{n};n\geq 1\right) $\ be a sequence of real numbers. From (\ref%
{x2}), it suffices to express $Z\left( n,s\right) $\ by considering the
binomial type sequence defined in (\ref{p5}),\ put after $b=a+x,\ c=a+ar+xr.$
To obtain (\ref{z4}), it suffices to choice $b=\left( r+1\right) c,\ c=x$\
in (\ref{z3}).
\end{proof}

\begin{example}
If $x_{n}=n,$\ the sequence given by (\ref{z4}) becomes%
\begin{equation*}
Z\left( n,s\right) =\frac{1}{\left( Df_{1}\left( 0\right) \right) ^{s}}%
\underset{j=0}{\overset{n}{\sum }}\binom{n}{j}\left( R+j\right)
^{n-j-1}\left( R+n\right) ^{j}x^{j}
\end{equation*}%
and if $x_{n}=n!,$\ the sequence given by (\ref{z4}) becomes%
\begin{equation*}
Z\left( n,s\right) =\frac{\left( R+n-1\right) !}{\left( Df_{1}\left(
0\right) \right) ^{s}}\underset{j=0}{\overset{n}{\sum }}\binom{n}{j}\frac{%
\left( R+n\right) ^{j}}{\left( R+j\right) !}x^{j}.
\end{equation*}
\end{example}

\subsection{Some\ applications\ of\ Theorem\ \protect\ref{3}}

Some particular cases of Theorem \ref{3} are given by the following
corollaries:

\begin{corollary}
Under the hypothesis of Theorem \ref{3} and $v\geq u$ the sequences%
\begin{eqnarray}
Y_{1}\left( n,k\right) &:&=\binom{n}{k}\frac{sk}{T}\underset{j=0}{\overset{T}%
{\sum }}\binom{T}{j}D_{z=\beta j+\lambda T}^{ju+vT}\left\{ e^{\alpha
z}f_{n-k}\left( z;a\right) \right\} x^{j}y^{T-j}\text{,}  \label{z8} \\
Y_{2}\left( n,k\right) &:&=\frac{n!}{k!}\frac{sk}{T}\underset{j=0}{\overset{T%
}{\sum }}\binom{T}{j}\frac{D_{z=\beta j+\lambda T}^{ju+vT}f_{vT+n-k}\left(
z;a\right) }{\left( vT+n-k\right) !}x^{j}y^{T-j}
\end{eqnarray}%
satisfy (\ref{h}).
\end{corollary}

\begin{proof}
Let $a_{1}=p,\ a_{2}=2q$\ and$\ a_{m}=0$\ for\ $m\geq 3.$\ To obtain (\ref%
{z8}) it suffices to replace the identity $B_{j,T}\left( p,2q,0,0,...\right)
=\binom{T}{j-T}p^{2T-j}q^{j-T}$ in (\ref{x5}) and in (\ref{x55}).
\end{proof}

Some particular cases of (\ref{z88}) are given by: \newline
For $r=u=0,\ v=1$\ or $u=v=r=0$ in (\ref{z8}), the sequences%
\begin{eqnarray}
Y_{1}\left( n,k\right) &:&=\binom{n}{k}\underset{j=0}{\overset{sk}{\sum }}%
\binom{sk}{j}D_{z=\beta j+\lambda k}^{sk}\left\{ e^{\alpha z}f_{n-k}\left(
z;a\right) \right\} x^{j}y^{sk-j},  \label{z88} \\
Y_{2}\left( n,k\right) &:&=\frac{n!}{k!}\underset{j=0}{\overset{sk}{\sum }}%
\binom{sk}{j}\frac{D_{z=\beta j+\lambda k}^{sk}f_{n+\left( s-1\right)
k}\left( z;a\right) }{\left( n+\left( s-1\right) k\right) !}x^{j}y^{sk-j}, \\
Y_{3}\left( n,k\right) &:&=\binom{n}{k}\underset{j=0}{\overset{sk}{\sum }}%
\binom{sk}{j}f_{n-k}\left( \beta j+\lambda k;a\right) x^{j}y^{sk-j}
\end{eqnarray}%
satisfy (\ref{h}), and if $s=1$ the last sequences give%
\begin{eqnarray}
&&B_{n,k}\left( \alpha y,...,m\left( \alpha yf_{m-1}\left( \lambda ;a\right)
+xD_{z=\beta }f_{m-1}\left( z;a\right) \right) ,...\right)  \label{p2} \\
&=&\binom{n}{k}\underset{j=0}{\overset{k}{\sum }}\binom{k}{j}%
D_{z=0}^{k}\left\{ e^{\alpha z}f_{n-k}\left( \left( \beta -\lambda \right)
j+\lambda k+z;a\right) \right\} x^{j}y^{k-j},  \notag
\end{eqnarray}%
\begin{eqnarray}
&&B_{n,k}\left( \left( x+y\right) Df_{1}\left( 0\right) ,...,yD_{z=\lambda
}f_{n}\left( z;a\right) +xD_{z=\beta }f_{n}\left( z;a\right) ,...\right)
\label{p3} \\
&=&\frac{1}{k!}\underset{j=0}{\overset{k}{\sum }}\binom{k}{j}%
D_{z=0}^{k}f_{n}\left( \left( \beta -\lambda \right) j+\lambda k+z;a\right)
x^{j}y^{k-j}  \notag
\end{eqnarray}%
\begin{eqnarray}
&&B_{n,k}\left( y+x,...,m\left( yf_{m-1}\left( \lambda ;a\right)
+xf_{m-1}\left( \beta ;a\right) \right) ,...\right)  \label{p4} \\
&=&\binom{n}{k}\underset{j=0}{\overset{k}{\sum }}\binom{k}{j}f_{n-k}\left(
\left( \beta -\lambda \right) j+\lambda k;a\right) x^{j}y^{k-j}  \notag
\end{eqnarray}

\begin{example}
For $f_{n}\left( x\right) =x^{n},\ a=0\ $and $u=0$\ or $1$\ in (\ref{z8})
the sequences%
\begin{eqnarray*}
Y_{1}\left( n,k\right) &=&\frac{sk}{T}\binom{n}{k}\underset{j=0}{\overset{T}{%
\sum }}\binom{T}{j}\left( \alpha T+\beta j\right) ^{n-k}x^{T-j}y^{j}\text{
and} \\
Y_{2}\left( n,k\right) &=&\frac{n!}{k!}\frac{sk}{T}\underset{j=0}{\overset{%
n-k}{\sum }}\binom{T}{j}\frac{\left( \alpha T+\beta j\right) ^{n-k-j}}{%
\left( n-k-j\right) !}x^{T-j}y^{j}
\end{eqnarray*}%
satisfy (\ref{h}), and for $f_{n}\left( x\right) =x^{n},\ a=0,\ u=0$\ or $%
1,\ a_{n}=n$\ in Theorem \ref{3} the sequences%
\begin{eqnarray*}
Y_{3}\left( n,k\right) &=&\frac{sk}{T}\binom{n}{k}\underset{j=0}{\overset{%
\infty }{\sum }}\left( \alpha T+\beta j\right) ^{n-k}\frac{\left( Tx\right)
^{j}}{j!}\text{ and} \\
Y_{4}\left( n,k\right) &=&\frac{sk}{T}\binom{n}{k}\underset{j=0}{\overset{n-k%
}{\sum }}\binom{n-k}{j}\left( \alpha T+\beta j\right) ^{n-k-j}\left(
Tx\right) ^{j}
\end{eqnarray*}%
satisfy (\ref{h}).
\end{example}

An interesting relation between Bell polynomials and Appell polynomials can
be viewed as a special case of Theorem \ref{3} and it is given by:

\begin{corollary}
For the sequence of polynomials $\left( A_{n}\left( x,y,z\right) \right) $
defined by%
\begin{equation*}
A_{n}\left( x,y,z\right) :=\underset{j=0}{\overset{n}{\sum }}\binom{n}{j}%
a_{j}\left( \left( x+y\right) j+x+y+z\right) \left( xj+yn+x+y+z\right)
^{n-1-j}
\end{equation*}%
we have%
\begin{eqnarray}
&&B_{n,k}\left( A_{0}\left( x,y,z\right) ,..,mA_{m-1}\left( x,y,z\right)
,...\right)  \label{z9} \\
&=&\underset{j=k}{\overset{n}{\sum }}\binom{n}{j}B_{j,k}\left(
a_{0},2a_{1},3a_{2},...\right) \left( jx+ny+kz\right) ^{n-j-1}\left( \left(
x+y\right) j+kz\right) ,  \notag
\end{eqnarray}%
and in particular when $x=y=0$\ we get%
\begin{equation}
B_{n,k}\left( A_{0}\left( z\right) ,...,mA_{m-1}\left( z\right) ,...\right) =%
\underset{j=k}{\overset{n}{\sum }}\binom{n}{j}B_{j,k}\left(
a_{0},...,ma_{m-1},...\right) \left( kz\right) ^{n-j},  \label{z10}
\end{equation}%
where $A_{n}\left( z\right) :=\underset{j=0}{\overset{n}{\sum }}\binom{n}{j}%
a_{j}z^{n-j}$\ is an Appell polynomial.
\end{corollary}

\begin{proof}
If suffices to use (\ref{x5}) with $x=1,\ r=v=0,\ s=u=1,$\ $f_{n}\left(
x;a\right) :=a\left( an+x\right) ^{n-1}$\ and replace $a_{n}$\ by $na_{n-1}.$
\end{proof}

\begin{example}
For any sequence $\left( \varphi _{n};n\geq 1\right) $\ for real numbers,
let $I_{n}$\ be the identity matrix of order $n\ $and $\left( A_{n}\right) $%
\ be the sequence of matrices defined by: $A_{0}:=1,\ A_{n}:=\left(
a_{ij}\right) $\ for $1\leq i,j\leq n$\ \ with $a_{ij}=\varphi _{j-i+1}$\ \
if \ $j\geq i,$\ $a_{i,i-1}=i-1$\ and $a_{ij}=0$\ otherwise. Then from $%
\left[ 7,\ p.\ 110\right] $\ and from (\ref{z10}) we get:%
\begin{equation*}
B_{n,k}\left( 1,..,m\det \left( A_{m-1}+xI_{m-1}\right) ,..\right) =\underset%
{j=k}{\overset{n}{\sum }}\binom{n}{j}B_{j,k}\left( 1,..,m\det
A_{m-1},..\right) \left( kx\right) ^{n-j}.
\end{equation*}
\end{example}

\subsection{Some\ applications\ of\ Theorem\ \protect\ref{4}}

A particular case of Theorem \ref{4} is given by the following corollary:

\begin{corollary}
Under the hypothesis of Theorem \ref{4} and $v\geq u$ the sequences%
\begin{eqnarray}
Z_{1}\left( n,s\right) &:&=\frac{1}{\gamma _{1}^{s}R}\underset{j=0}{\overset{%
R}{\sum }}\binom{R}{j}D_{z=\beta j+\lambda R}^{ju+vR}\left\{ e^{\alpha
z}f_{n}\left( z;a\right) \right\} x^{j}y^{R-j}\text{ and}  \label{z12} \\
Z_{2}\left( n,s\right) &:&=\frac{n!}{\gamma _{2}^{s}R}\underset{j=0}{\overset%
{R}{\sum }}\binom{R}{j}\frac{D_{z=\beta j+\lambda R}^{ju+vR}f_{n+\left(
u+v\right) R}\left( z;a\right) }{\left( n+\left( u+v\right) R\right) !}%
x^{j}y^{R-j}
\end{eqnarray}%
satisfy (\ref{b}), where $\gamma _{1}:=\alpha ^{v}\left( x\alpha
^{u}+y\right) $\ and $\gamma _{2}:=\left\{ 
\begin{array}{cc}
\left( x+y\right) \left( Df_{1}\left( 0\right) \right) ^{v} & \text{if }%
u\geq 1 \\ 
y\left( Df_{1}\left( 0\right) \right) ^{v} & \text{if }u=0%
\end{array}%
\right. .$
\end{corollary}

\begin{proof}
It suffices to put in Theorem \ref{4} $x=1;\ a_{1}=p,\ a_{2}=2q$ and$\
a_{m}=0$\ for $m\geq 3$\ and use the identity $B_{j,R}\left(
p,2q,0,0,...\right) =\binom{R}{j-R}p^{2R-j}q^{j-R}$.
\end{proof}

\begin{example}
For $f_{n}\left( x\right) =x^{n},\ a=0$\ and $u=0$\ or $1$\ the sequence
given by (\ref{z12}) becomes%
\begin{eqnarray*}
Z_{1}\left( n,s\right) &=&\left( x+y\right) ^{-s}\frac{1}{R}\underset{j=0}{%
\overset{R}{\sum }}\binom{R}{j}\left( \beta j+\alpha R\right)
^{n}x^{j}y^{R-j}\text{ and} \\
Z_{2}\left( n,s\right) &=&\frac{n!}{R}\underset{j=0}{\overset{n}{\sum }}%
\frac{1}{\left( n-j\right) !}\binom{R}{j}\left( \beta j+\alpha R\right)
^{n-j}x^{j}
\end{eqnarray*}%
satisfy (\ref{b}), and for $f_{n}\left( x\right) =x^{n},\ a=0,\ x=1,\ u=0$\
or $1,\ a_{n}=nt^{n-1}$\ the sequence given by (\ref{x7})becomes%
\begin{eqnarray*}
Z_{3}\left( n,s\right) &=&\frac{1}{R}\underset{j=0}{\overset{n}{\sum }}%
\binom{n}{j}\left( \beta j+\alpha R\right) ^{n-j}R^{j}t^{j}\text{ and} \\
Z_{4}\left( n,s\right) &=&\frac{\exp \left( -sx\right) }{R!}\underset{j=0}{%
\overset{\infty }{\sum }}\left( \beta j+\alpha R\right) ^{n}\frac{\left(
Rx\right) ^{j}}{j!}.
\end{eqnarray*}
\end{example}

\section{Proof of the main results}

\begin{proof}[Proof of Lemma 1]
Let $F\left( t;a\right) ^{x}:=1+\overset{\infty }{\underset{n=1}{\sum }}%
f_{n}\left( x;a\right) \frac{t^{n}}{n!}$. Now, because $f_{n}\left(
x;a\right) $ is a polynomial of degree $n,$ then the proof follows from the
following expansions:%
\begin{equation*}
\overset{\infty }{\underset{n=0}{\sum }}D_{z=0}^{k}\left( e^{\alpha
z}f_{n}\left( kx+z;a\right) \right) \frac{t^{n+k}}{n!}=\overset{\infty }{%
\underset{n=k}{\sum }}\frac{n!}{\left( n-k\right) !}D_{z=0}^{k}\left(
e^{\alpha z}f_{n-k}\left( kx+z;a\right) \right) \frac{t^{n}}{n!}
\end{equation*}%
\newline
and%
\begin{eqnarray*}
\overset{\infty }{\underset{n=0}{\sum }}D_{z=0}^{k}\left( e^{\alpha
z}f_{n}\left( kx+z;a\right) \right) \frac{t^{n+k}}{n!} &=&D_{z=0}^{k}\left(
e^{\alpha z}\overset{\infty }{\underset{n=0}{\sum }}f_{n}\left(
kx+z;a\right) \frac{t^{n+k}}{n!}\right) \\
&=&D_{z=0}^{k}\left( e^{\alpha z}F\left( t;a\right) ^{kx+z}t^{k}\right) \\
&=&t^{k}F\left( t;a\right) ^{kx}D_{z=0}^{k}\left( e^{\left( \alpha +\ln
F\left( t;a\right) \right) z}\right) \\
&=&t^{k}F\left( t;a\right) ^{kx}\left( \alpha +\ln F\left( t;a\right)
\right) ^{k} \\
&=&t^{k}e^{-\alpha kx}\left( \left( e^{\alpha }F\left( t;a\right) \right)
^{x}\ln \left( e^{\alpha }F\left( t;a\right) \right) \right) ^{k} \\
&=&t^{k}e^{-\alpha kx}\left( D_{x}\left( e^{\alpha }F\left( t;a\right)
\right) ^{x}\right) ^{k} \\
&=&t^{k}\left( \overset{\infty }{\underset{m=0}{\sum }}e^{-\alpha
x}D_{x}\left( e^{\alpha x}f_{m}\left( x;a\right) \right) \frac{t^{m}}{m!}%
\right) ^{k} \\
&=&\left( \overset{\infty }{\underset{m=1}{\sum }}me^{-\alpha x}D_{x}\left(
e^{\alpha x}f_{m-1}\left( x;a\right) \right) \frac{t^{m}}{m!}\right) ^{k} \\
&=&k!\overset{\infty }{\underset{n=k}{\sum }}B_{n,k}\left( \alpha
,...,me^{-\alpha x}D_{x}\left( e^{\alpha x}f_{m-1}\left( x;a\right)
,...\right) \right) \frac{t^{n}}{n!}.
\end{eqnarray*}%
\newline
Then by identification we get%
\begin{equation*}
B_{n,k}\left( \alpha ,2e^{-\alpha x}D_{x}\left( e^{\alpha x}f_{1}\left(
x;a\right) \right) ,...,mD_{x}\left( e^{\alpha x}f_{m-1}\left( x;a\right)
\right) \right) =\binom{n}{k}D_{z=0}^{k}\left( e^{\alpha z}f_{n-k}\left(
kx+z;a\right) \right) .
\end{equation*}%
To obtain (\ref{alpha}), it suffices to remark that for $m\geq 1$ we have%
\begin{equation*}
e^{-\alpha x}D_{x}\left( e^{\alpha x}f_{m}\left( x;a\right) \right)
=D_{z=0}\left( e^{\alpha z}f_{m}\left( x+z;a\right) \right) .
\end{equation*}%
The identity (\ref{alpha}) can be written when $\alpha =0$ as%
\begin{equation*}
B_{n,k}\left( 0,...,mD_{z=0}f_{m-1}\left( kx+z;a\right) ,...\right) =\binom{n%
}{k}D_{z=0}^{k}f_{n-k}\left( kx+z;a\right)
\end{equation*}%
which is equivalent to%
\begin{equation*}
B_{n,k}\left( D_{x}f_{1}\left( x;a\right) ,...,D_{x}f_{m}\left( x;a\right)
,...\right) =\frac{1}{k!}D_{z=0}^{k}f_{n}\left( kx+z;a\right) .
\end{equation*}
\end{proof}

\begin{proof}[Proof of Theorem \protect\ref{1}]
When we replace $x_{n}$ by $\alpha x_{n}$ in (\ref{w1}) and we use the
well-known identities%
\begin{eqnarray}
B_{n,k}\left( \alpha x_{1},\alpha x_{2},\alpha x_{3}...\right) &=&\alpha
^{k}B_{n,k}\left( x_{1},x_{2},x_{3},...\right) \text{ \ } \\
\text{and \ }B_{n,k}\left( \alpha x_{1},\alpha ^{2}x_{2},\alpha
^{3}x_{3}...\right) &=&\alpha ^{n}B_{n,k}\left( x_{1},x_{2},x_{3},...\right)
,
\end{eqnarray}%
it results that the identity (\ref{w1}) remains true for $x_{1}\neq 1$. Then
for $r+s\geq 1$ and for the choice $x_{n}=nD_{z=0}\left( e^{\alpha
z}f_{n-1}\left( x+z;a\right) \right) $\ in (\ref{w1}), the identity (\ref%
{alpha}) proves that the sequence $\left( Y\left( n,k\right) \right) $ given
by (\ref{w2}) becomes%
\begin{eqnarray*}
Y\left( n,k\right) &=&\binom{n}{k}\frac{sk}{T}\binom{T+n-k}{T}%
^{-1}B_{T+n-k,\ T}\left( \alpha ,...,mD_{z=0}\left( e^{\alpha
z}f_{m-1}\left( x+z;a\right) \right) ,...\right) \\
&=&\binom{n}{k}\frac{sk}{T}D_{z=0}^{T}\left\{ e^{\alpha z}f_{n-k}\left(
Tx+z;a\right) \right\} .
\end{eqnarray*}%
For the particular case $\alpha =0,$ if we take $x_{n}=D_{z=0}\left(
f_{n}\left( x+y;a\right) \right) $ in (\ref{w1}) and we use the identity (%
\ref{alpha1}), the sequence $\left( Y\left( n,k\right) \right) $ given by (%
\ref{w2}) becomes%
\begin{eqnarray*}
Y\left( n,k\right) &=&\binom{n}{k}\frac{sk}{T}\binom{T+n-k}{T}%
^{-1}B_{T+n-k,\ T}\left( D_{z=0}\left( e^{\alpha z}f_{1}\left( x+z;a\right)
\right) ,...\right) \\
&=&\frac{n!}{k!}\frac{sk}{T}\frac{1}{\left( T+n-k\right) !}%
D_{z=0}^{T}f_{T+n-k}\left( Tx+z;a\right) .
\end{eqnarray*}%
Note that for the case $r=s=0,$ the sequence $\left( Y\left( n,k\right)
\right) $ given by (\ref{x11}) is not defined. We put in this case $Y\left(
n,k\right) =\binom{n}{k}f_{n-k}\left( x;a\right) $. Proposition 1 in $\left[
4\right] $ proves that this sequence satisfies (\ref{h}).
\end{proof}

\begin{proof}[Proof of Theorem \protect\ref{2}]
\textbf{Case }$\alpha \neq 0:$ Let $x_{n}=\frac{n}{\alpha }D_{z=0}\left\{
e^{\alpha z}f_{n-1}\left( x+z;a\right) \right\} .$ We have $x_{1}=1$ and
then by using the identity (\ref{alpha}), the sequence $\left( Z\left(
n,s\right) \right) $ given by (\ref{w4}) becomes%
\begin{eqnarray*}
Z\left( n,s\right) &=&\frac{1}{\alpha ^{R}}\frac{1}{R}\binom{R+n}{R}%
^{-1}B_{R+n,\ R}\left( \alpha ,...,mD_{z=0}\left( e^{\alpha z}f_{m-1}\left(
x+z;a\right) \right) ,...\right) \\
&=&\frac{1}{\alpha ^{R}}\frac{1}{R}D_{z=0}^{R}\left( e^{\alpha z}f_{n}\left(
Rx+z;a\right) \right) .
\end{eqnarray*}%
\textbf{Case }$\alpha =0:$ Let $x_{n}=D_{x}f_{n}\left( x;a\right)
/D_{x}f_{1}\left( x;a\right) .$\ Now because $D_{x}f_{1}\left( x\right)
=xDf_{1}\left( 0\right) $ we get $D_{x}f_{1}\left( x;a\right) =D_{x}\left( 
\frac{x}{x+a}f_{1}\left( x+a\right) \right) =Df_{1}\left( 0\right) \neq 0.$
We have $x_{1}=1\ $and then by using the identity (\ref{alpha1}), the
sequence $\left( Z\left( n,s\right) \right) $ given by (\ref{w4}) becomes%
\begin{eqnarray*}
Z\left( n,s\right) &=&\frac{1}{R}\binom{R+n}{R}^{-1}\frac{B_{R+n,\ R}\left(
D_{z=0}\left( f_{1}\left( x+z;a\right) \right) ,...,D_{z=0}\left(
f_{m}\left( x+z;a\right) \right) ,...\right) }{\left( D_{x}f_{1}\left(
x;a\right) \right) ^{R}} \\
&=&\frac{1}{\left( Df_{1}\left( 0\right) \right) ^{R}}\frac{1}{R}\frac{n!}{%
\left( R+n\right) !}D_{z=0}^{R}f_{R+n}\left( kx+z;a\right) .
\end{eqnarray*}%
Note that if $Z\left( n,s\right) $ satisfies (\ref{b}) then $\lambda
^{n}Z\left( n,s\right) $ satisfies (\ref{b}).
\end{proof}

\begin{proof}[Proof of Theorem \protect\ref{3}]
Let $\left( a_{n};n\geq 1\right) $\ be a real sequences; $u,v$ be a natural
numbers and we put $F\left( t\right) ^{x}:=1+\overset{\infty }{\underset{n=1}%
{\sum }}f_{n}\left( x\right) \frac{t^{n}}{n!}\ $with $f_{0}\left( x\right)
:=1.$\newline
For $F\left( n,k\right) :=\underset{j\geq k}{\sum }B_{j,k}\left(
a_{1},a_{2},\ldots \right) D_{z=0}^{ju+vk}\left\{ e^{\alpha z}f_{n}\left(
\beta j+\lambda k+z;a\right) \right\} \frac{x^{j}}{j!}$\ we have%
\begin{eqnarray*}
\underset{n\geq 0}{\sum }F\left( n,k\right) \frac{t^{n+k}}{n!} &=&t^{k}%
\underset{n\geq 0}{\sum }\frac{t^{n}}{n!}\left( \underset{j\geq k}{\sum }%
B_{j,k}\left( a_{1},a_{2},\ldots \right) D_{z=0}^{ju+vk}\left\{ e^{\alpha
z}f_{n}\left( \beta j+\lambda k+z;a\right) \right\} \frac{x^{j}}{j!}\right)
\\
&=&t^{k}\underset{j\geq k}{\sum }B_{j,k}\left( a_{1},a_{2},\ldots \right) 
\dfrac{x^{j}}{j!}\left( \underset{n\geq 0}{\sum }D_{z=0}^{ju+vk}\left\{
e^{\alpha z}f_{n}\left( \beta j+\lambda k+z;a\right) \right\} \frac{t^{n}}{n!%
}\right) \\
&=&t^{k}\underset{j\geq k}{\sum }B_{j,k}\left( a_{1},a_{2},\ldots \right) 
\frac{x^{j}}{j!}D_{z=0}^{ju+vk}\left\{ e^{\alpha z}\underset{n\geq 0}{\sum }%
f_{n}\left( \beta j+\lambda k+z;a\right) \dfrac{t^{n}}{n!}\right\} \\
&=&t^{k}\underset{j\geq k}{\sum }B_{j,k}\left( a_{1},a_{2},\ldots \right) 
\frac{x^{j}}{j!}D_{z=0}^{ju+vk}\left\{ e^{\alpha z}\left( F\left( t;a\right)
\right) ^{\beta j+\lambda k+z}\right\} \\
&=&t^{k}\underset{j\geq k}{\sum }B_{j,k}\left( a_{1},a_{2},\ldots \right) 
\frac{x^{j}}{j!}\left( F\left( t;a\right) \right) ^{\beta j+\lambda
k}D_{z=0}^{ju+vk}\left\{ e^{\alpha z}F\left( t;a\right) ^{z}\right\} \\
&=&t^{k}\underset{j\geq k}{\sum }B_{j,k}\left( a_{1},a_{2},\ldots \right) 
\frac{x^{j}}{j!}\left( F\left( t;a\right) \right) ^{\beta j+\lambda k}\left(
\alpha +\ln F\left( t;a\right) \right) ^{ju+vk} \\
&=&t^{k}\underset{j\geq k}{\sum }B_{j,k}\left( u_{1},u_{2},\ldots \right) 
\frac{x^{j}}{j!} \\
&=&\frac{t^{k}}{k!}\left( \underset{m\geq 1}{\sum }u_{m}\frac{x^{m}}{m!}%
\right) ^{k} \\
&=&\frac{t^{k}}{k!}\left( \underset{m\geq 1}{\sum }a_{m}F\left( t;a\right)
^{\beta m+\lambda }\left( D_{z=0}\left( e^{\alpha z}\left( F\left(
t;a\right) \right) ^{z}\right) \right) ^{mu+v}\frac{x^{m}}{m!}\right) ^{k} \\
&=&\frac{t^{k}}{k!}\left( \underset{m\geq 1}{\sum }a_{m}D_{z=0}^{mu+v}\left(
e^{\alpha z}F\left( t;a\right) ^{\beta m+\lambda +z}\right) \frac{x^{m}}{m!}%
\right) ^{k} \\
&=&\frac{t^{k}}{k!}\left( \underset{m\geq 1}{\sum }a_{m}\underset{j\geq 0}{%
\sum }D_{z=0}^{mu+v}\left( e^{\alpha z}f_{j}\left( \beta m+\lambda
+z;a\right) \frac{t^{j}}{j!}\right) \frac{x^{m}}{m!}\right) ^{k} \\
&=&\frac{t^{k}}{k!}\left( \underset{j\geq 0}{\sum }\frac{t^{j}}{j!}\underset{%
m\geq 1}{\sum }a_{m}D_{z=0}^{mu+v}\left( e^{\alpha z}f_{j}\left( \beta
m+\lambda +z;a\right) \right) \frac{x^{m}}{m!}\right) ^{k} \\
&=&\frac{t^{k}}{k!}\left( \underset{j\geq 0}{\sum }v_{j}\frac{t^{j}}{j!}%
\right) ^{k} \\
&=&\frac{1}{k!}\left( \underset{j\geq 1}{\sum }jv_{j-1}\frac{t^{j}}{j!}%
\right) ^{k} \\
&=&\underset{n\geq k}{\sum }B_{n,k}\left( v_{0},...,jv_{j-1},...\right) 
\frac{t^{n}}{n!},\ \ \ 
\end{eqnarray*}%
where $u_{m}=a_{m}\left( F\left( t;a\right) \right) ^{\beta m+\lambda
}\left( \alpha +\ln F\left( t;a\right) \right) ^{mu+v}=a_{m}F\left(
t;a\right) ^{\beta m+\lambda }\left( D_{z=0}\left( e^{\alpha z}\left(
F\left( t;a\right) \right) ^{z}\right) \right) ^{mu+v}\ \ \left( m\geq
1\right) $\ \ and \ $v_{m}=\underset{m\geq 1}{\sum }a_{m}D_{z=0}^{mu+v}%
\left( e^{\alpha z}f_{j}\left( \beta m+\lambda +z;a\right) \right) \frac{%
x^{m}}{m!}=F\left( m,1\right) \ \ \left( m\geq 0\right) .$ Then, we obtain%
\begin{equation*}
B_{n,k}\left( F\left( 0,1\right) ,\ldots ,jF\left( j-1,1\right) ,\ldots
\right) =k!\binom{n}{k}F\left( n-k,k\right)
\end{equation*}%
and this means that the sequence%
\begin{equation}
X\left( n,k\right) :=k!\binom{n}{k}\underset{j\geq k}{\sum }B_{j,k}\left(
a_{1},a_{2},\ldots \right) D_{z=0}^{ju+vk}\left\{ e^{\alpha z}f_{n-k}\left(
\beta j+\lambda k+z;a\right) \right\} \frac{x^{j}}{j!}  \label{y5}
\end{equation}%
satisfies (\ref{h}). To obtain (\ref{x55}), it suffices to take $%
x_{n}=X\left( n,1\right) $ in (\ref{w1}), where $\left( X\left( n,k\right)
\right) $ is the sequence given by (\ref{y5}). Indeed, the sequence given by
(\ref{w2}) becomes%
\begin{eqnarray*}
Y\left( n,k\right) &=&\binom{n}{k}\frac{sk}{T}\binom{T+n-k}{T}%
^{-1}B_{T+n-k,\ T}\left( X\left( 1,1\right) ,X\left( 2,1\right) ,\ldots
\right) \\
&=&\binom{n}{k}\frac{sk}{T}\binom{T+n-k}{T}^{-1}X\left( T+n-k,\ T\right) \\
&=&\binom{n}{k}\frac{sk}{T}T!\underset{j\geq T}{\sum }B_{j,T}\left(
a_{1},a_{2},\ldots \right) D_{z=0}^{ju+vT}\left\{ e^{\alpha z}f_{n-k}\left(
\beta j+\lambda T+z;a\right) \right\} \frac{x^{j}}{j!}.
\end{eqnarray*}%
\newline
For the particular case $\alpha =0$ and $u\geq 1$ we remark that $X\left(
n,1\right) =0$ for $u+v>n-1$ and then the identity (\ref{y5}) becomes 
\begin{equation}
B_{n,k}\left( \overset{u+v}{\overbrace{0,...,0}},X\left( u+v+1,1\right)
,X\left( u+v+2,1\right) ,\ldots \right) =X\left( n,k\right)  \label{yyy}
\end{equation}%
or, equivalently, by using the well-known identity%
\begin{equation*}
B_{n,k}\left( 0,...,0,a_{r+1},a_{r+2},...\right) =\frac{n!}{\left(
n-rk\right) !}B_{n-rk,k}\left( \frac{a_{1+r}}{\left( 1+r\right) !},\ldots ,%
\frac{i!a_{i+r}}{\left( i+r\right) !},\ldots \right) .
\end{equation*}%
Setting $X^{\ast }\left( n,k\right) :=n!X\left( n+\left( u+v\right)
k,k\right) /\left( n+\left( u+v\right) k\right) !,$ the identity (\ref{yyy})
becomes%
\begin{equation}
B_{n,\ k}\left( X^{\ast }\left( 1,1\right) ,X^{\ast }\left( 2,1\right)
,\ldots \right) =X^{\ast }\left( n,k\right)  \notag
\end{equation}%
To obtain (\ref{x5}) it suffices to take $x_{n}=X^{\ast }\left( n,1\right) $
in (\ref{w1}). Indeed%
\begin{eqnarray*}
Y\left( n,k\right) &=&\binom{n}{k}\frac{sk}{T}\binom{T+n-k}{T}%
^{-1}B_{T+n-k,\ T}\left( X^{\ast }\left( 1,1\right) ,X^{\ast }\left(
2,1\right) ,\ldots \right) \\
&=&\binom{n}{k}\frac{sk}{T}\binom{T+n-k}{T}^{-1}X^{\ast }\left( T+n-k,\
T\right) \\
&=&\frac{n!}{k!}\frac{sk}{T}T!\underset{j\geq T}{\sum }B_{j,T}\left(
a_{1},a_{2},\ldots \right) D_{z=0}^{ju+vT}\left\{ e^{\alpha z}\frac{%
f_{\left( u+v\right) T+n-k}\left( \beta j+\lambda T+z;a\right) }{\left(
\left( u+v\right) T+n-k\right) !}\right\} \frac{x^{j}}{j!}.
\end{eqnarray*}%
\newline
\end{proof}

\begin{proof}[Proof of Theorem \protect\ref{4}]
For $\alpha \neq 0$ it suffices to put in (\ref{w3}) $x_{n}=X\left(
n,1\right) /X\left( 1,1\right) ,$\ where $\left( X\left( n,k\right) \right) $%
\ is the sequence given by (\ref{y5}). We have 
\begin{equation*}
x_{1}=X\left( 1,1\right) =\underset{j\geq 1}{\sum }a_{j}\alpha ^{ju+v}\frac{%
x^{j}}{j!}=\alpha ^{v}\varphi \left( x\alpha ^{u}\right) \text{ \ with }%
\varphi \left( x\right) :=\underset{j\geq 1}{\sum }a_{j}\frac{x^{j}}{j!}.
\end{equation*}%
The sequence $\left( Z\left( n,s\right) \right) $ given by (\ref{w4}) becomes%
\begin{eqnarray*}
Z\left( n,s\right) &=&\frac{1}{\left( X\left( 1,1\right) \right) ^{R}}\frac{1%
}{R}\binom{R+n}{R}^{-1}X\left( R+n,R\right) \\
&=&\frac{1}{\left( \alpha ^{v}\varphi \left( x\alpha ^{u}\right) \right) ^{R}%
}\frac{R!}{R}\underset{j\geq R}{\sum }B_{j,R}\left( a_{1},a_{2},\ldots
\right) D_{z=0}^{ju+vR}\left\{ e^{\alpha z}f_{n}\left( \beta j+\lambda
R+z;a\right) \right\} \frac{x^{j}}{j!}.
\end{eqnarray*}%
\newline
For $\alpha =0$ it suffices to put in (\ref{w3}) $x_{n}=X^{\ast }\left(
n,1\right) /X^{\ast }\left( 1,1\right) ,$\ where $\left( X^{\ast }\left(
n,k\right) \right) $\ is the sequence given by (\ref{y55}). We have%
\begin{equation*}
X^{\ast }\left( 1,1\right) =\frac{X\left( u+v+1,1\right) }{\left(
u+v+1\right) !}=\left\{ 
\begin{array}{cc}
a_{1}x\left( Df_{1}\left( 0\right) \right) ^{u+v} & \text{if }u\geq 1 \\ 
\left( Df_{1}\left( 0\right) \right) ^{v}\varphi \left( x\right) & \text{if }%
u=0%
\end{array}%
\right.
\end{equation*}%
The sequence $\left( Z\left( n,s\right) \right) $ given by (\ref{w4}) becomes%
\begin{eqnarray*}
Z\left( n,s\right) &=&\frac{1}{\left( X^{\ast }\left( 1,1\right) \right) ^{R}%
}\frac{1}{R}\binom{R+n}{R}^{-1}X^{\ast }\left( R+n,R\right) \\
&=&\frac{n!}{\left( X^{\ast }\left( 1,1\right) \right) ^{R}}\frac{R!}{R}%
\underset{j\geq R}{\sum }B_{j,R}\left( a_{1},a_{2},\ldots \right)
D_{z=0}^{ju+vR}\left\{ \frac{f_{n+\left( u+v\right) R}\left( \beta j+\lambda
R+z;a\right) }{\left( n+\left( u+v\right) R\right) !}\right\} \frac{x^{j}}{j!%
}.
\end{eqnarray*}%
\newline
\end{proof}

\begin{acknowledgement}
The author thanks the Professor Hacene Belbachir for his careful reading and
valuable suggestions.
\end{acknowledgement}

\end{document}